\documentclass[11pt]{article}

\usepackage[utf8]{inputenc}
\usepackage[T1]{fontenc}
\usepackage[american]{babel}
\usepackage{etex}

\usepackage[margin=3cm]{geometry}
\usepackage{tikz}
\usetikzlibrary{patterns}
\usepackage{csquotes}
\usepackage{amsmath,amsfonts,amstext,amssymb,amsbsy,amsopn,eucal,dsfont,mathtools}
\mathtoolsset{mathic}
\usepackage[amsmath,thmmarks]{ntheorem}

\usepackage{graphicx}
\usepackage[shortlabels]{enumitem}
\usepackage[vlined,ruled]{algorithm2e}
\usepackage{xspace}

\usepackage{hyperref}
\usepackage[capitalize]{cleveref}

\usepackage{microtype}

\newlist{compactenum}{enumerate}{4}
\setlist[compactenum]{label=\arabic*.,itemsep=0pt,topsep=0pt,parsep=0pt,leftmargin=*}
\newlist{compactitem}{itemize}{4}
\setlist[compactitem]{label=\textbullet,itemsep=0pt,topsep=0pt,parsep=0pt}

\crefformat{enumi}{#2#1#3} \crefrangeformat{enumi}{#3#1#4--#5#2#6}
\crefmultiformat{enumi}{#2#1#3}{ and~#2#1#3}{, #2#1#3}{ and~#2#1#3}
\crefrangemultiformat{enumi}{#3#1#4--#5#2#6}{
and~#3#1#4--#5#2#6}{,#3#1#4--#5#2#6}{ and~#3#1#4--#5#2#6}
\crefalias{enumii}{enumi} \crefalias{enumiii}{enumi}
\crefalias{compactenumi}{enumi} \crefalias{compactenumii}{enumi}
\crefalias{compactenumiii}{enumi} \crefformat{equation}{(#2#1#3)}
\crefrangeformat{equation}{(#3#1#4)--(#5#2#6)}
\crefmultiformat{equation}{(#2#1#3)}{ and~(#2#1#3)}{, (#2#1#3)}{ and~(#2#1#3)}
\crefrangemultiformat{equation}{(#3#1#4)--(#5#2#6)}{ and~(#3#1#4)--(#5#2#6)}{,
(#3#1#4)--(#5#2#6)}{ and~(#3#1#4)--(#5#2#6)}

\newtheorem{Theorem}{Theorem}
\newtheorem{Lemma}[Theorem]{Lemma}

\newtheorem{Corollary}[Theorem]{Corollary}

\newtheorem{Definition}[Theorem]{Definition}

\newtheorem{Problem-}[Theorem]{Problem}

\newtheorem{Example}[Theorem]{Example}

\theoremstyle{proof} \theorembodyfont{\normalfont}
\theoremsymbol{\ensuremath{\square}}
\newtheorem*{Proof}{Proof}
\author{Christoph Helmberg\thanks{Fakult\"at f\"ur Mathematik, Technische
    Universit\"at Chemnitz, D-09107 Chemnitz,
    Germany. helmberg@mathematik.tu-chemnitz.de} \and Vilmar
  Trevisan\thanks{Instituto de Matem\'atica, Universidade Federal do
    Rio Grande do Sul, CEP 91509-900, Porto Alegre, RS, Brazil. trevisan@mat.ufrgs.br}%
}%
\title{Spectral threshold dominance, Brouwer's conjecture and maximality of Laplacian energy}
\date{\today}

\begin{document}

\maketitle
\begin{abstract}
The Laplacian energy of a graph is the sum of the distances of the
eigenvalues of the Laplacian matrix of the graph to the graph's average
degree. The maximum Laplacian energy over all graphs on $n$ nodes
and $m$ edges is conjectured to be attained for threshold graphs.
We prove the conjecture to hold for graphs with the property that for
each $k$ there is a threshold graph on the same number of nodes
and edges whose sum of the $k$ largest Laplacian eigenvalues exceeds that
of the $k$ largest Laplacian eigenvalues of the graph. We call such graphs
spectrally threshold dominated. These graphs include split graphs
and cographs and spectral threshold dominance is preserved by
disjoint unions and taking complements. We conjecture that all graphs
are spectrally threshold dominated. This conjecture turns out to be
equivalent to Brouwer's conjecture concerning a bound on the
sum of the $k$ largest Laplacian eigenvalues.


\end{abstract}

\section{Introduction}
Let $G=(N,E)$ be a simple graph with node set $N=\{1,\dots,n\}$ and
edge set $E\subseteq\{\{i,j\}:i,j\in N,i\neq j\}$. For brevity, we
will usually write $ij$ instead of $\{i,j\}$ for edges and put
$m=|E|$. It will be convenient to assume that the nodes are numbered
so that their degrees $d_i=|\{j: ij\in E\}|$ are sorted
nonincreasingly.  Let $e_i$ denote the $i$-th column of the $n\times
n$ identity matrix $I_n$ and define the positive semidefinite matrices
$E_{ij}:=(e_i-e_j)(e_i-e_j)^T$, then the Laplacian matrix of $G$ is
defined to be $L(G)=\sum_{ij\in E}E_{ij}$. $L(G)$ may also be written
as $L(G) = D(G) - A(G)$, where $D(G)$ is the diagonal degree matrix
and $A(G)$ is the adjacency matrix of $G$. If $G$ is clear from the
context, we drop the argument and simply write $L$. The Laplacian is a
positive semidefinite matrix with a trivial eigenvalue $0$ and the
vector of all ones $\mathbf{1}$ as associated eigenvector. In this
paper we denote the eigenvalues of $L$ in nonincreasing order by
$\lambda_1(L)\ge\dots\ge\lambda_{n-1}(L)\ge \lambda_n(L)=0$.  As the
trace of $L$ is $2m$ there holds $\sum_{i=1}^n\lambda_i(L)=2m$ and for
$m>0$ at least one eigenvalue has value greater than the average
degree $2m/n$. The eigenvalues of $L$ are the Laplacian eigenvalues of $G$ and we may write $\lambda_i(L) = \lambda_i(G)$.

The Laplacian energy of a graph $G$ is defined to be
\begin{displaymath}
LE(G):=\sum_{i=1}^n\left|\lambda_i(G)-\frac{2m}n\right|
\end{displaymath}
and \cite{Aim} raised the question which graphs on $n$ nodes maximize
this value.

For $i=1,\dots,n$ the conjugate degree $d_i^*(G)=|\{i\colon d_i\ge
i\}|$ gives the number of nodes of $G$ of degree at least $i$.  Each
degree sequence satisfying $d_i^*=d_i+1$ for $i=1,\dots,f$ with
\emph{trace} $f=\max\{i\colon d_i\ge i\}$ uniquely defines a graph and
these graphs form the so called threshold graphs
\cite{MerrisRoby2005}. We notice that the degree sequence $d$ of
threshold graphs is fully specified once the conjugate degrees $d_i^*$
are given for $i\in[f]$, which is easily seen by looking at a Ferrers
diagram of $d$ (see next section for definitions).  There the part
strictly below the diagonal boxes is the transpose of the part above
and including the diagonal.  The right hand side of Figure
\ref{fig2ex} is the Ferrers diagram of the (threshold) graph
corresponding to the degree sequence $d=(4,4,3,3,2,0,0,0)$. The dual
degrees give the number of boxes of the columns.  Threshold graphs may
also be characterized combinatorially by starting with the empty graph
and iteratively adding a node that is isolated or fully connected to
all previous nodes.

In our context, a central property of threshold graphs $T$ is that the conjugate
degrees are exactly the eigenvalues of their Laplacian matrix,
$\lambda_i(T)=d_i^*(T)$ for $i=1,\dots,n$ \cite{Merris94}. It has been
conjectured that among all connected graphs
on $n$ nodes the threshold graph called pineapple with trace
$\lfloor \frac{2n}{3} \rfloor$ maximizes the Laplacian energy (see
\cite{Vin2013}). Among connected threshold graphs the pineapple is
indeed the maximizer; for general threshold graphs on $n$ nodes the
clique of size $\lfloor \frac{2n+1}{3} \rfloor +1 $ together with
$\lfloor \frac{n-3}{3} \rfloor$ isolated vertices is
a threshold graph maximizing Laplacian energy \cite{Helm15} and we
conjecture that this graph has maximum Laplacian energy among all
graphs on $n$ nodes.

In this paper we prove that the general conjecture holds for graphs that are
spectrally dominated by threshold graphs in the following sense.
\begin{Definition} A graph $G$ on $n$ nodes with $m$
  edges is  \emph {spectrally threshold dominated} if for each $k\in\{1,\dots,n\}$ there is a threshold graph $T_k$ having the  same number of nodes and edges satisfying $\sum_{i=1}^kd_i^*(T_k) =\sum_{i=1}^k\lambda_i(T_k) \ge \sum_{i=1}^k\lambda_i(G)$.
\end{Definition}
This definition was in part motivated by the Grone-Merris conjecture,
proved by Bai \cite{Bai11} -- from here on called the Grone-Merris-Bai
Theorem --  which states that for any graph $G$ on $n$ vertices with
degree sequence $d_1 \geq \ldots \geq d_n$ and for any $k \in \{1,\ldots ,n\}$,
\begin{equation}\label{Bai}
\sum_{i=1}^k\lambda_i(G) \leq \sum_{i=1}^kd_i^*(G).
\end{equation}
Note that equality holds in (\ref{Bai}) for threshold graphs.

Our main result (proved in \cref{S:specdom}) is the following.
\begin{Theorem}\label{Th:specdom}
  For each spectrally threshold dominated graph $G$ there exists a threshold
  graph with the same number of nodes and edges whose Laplacian
  energy is at least as large as that of $G$.
\end{Theorem}
We conjecture that all graphs are spectrally threshold dominated, in
which case the maximum Laplacian energy would be attained by
threshold graphs for any given number of nodes and edges. We prove that this class goes
well beyond threshold graphs (definitions of the graph classes will be
given along with the proofs in \cref{S:split} and \cref{S:preserve}).
\begin{Theorem}\label{Th:split}
  Split graphs are spectrally threshold dominated.
\end{Theorem}
\begin{Theorem}\label{Th:preserve}
  Disjoint unions and complements of spectrally threshold dominated
  graphs are spectrally threshold dominated.
\end{Theorem}
This has the following immediate consequence.
\begin{Corollary}\label{Th:cographs}
  Cographs are spectrally threshold dominated.
\end{Corollary}

The search for further examples of graph classes whose sum of the $k$
largest Laplacian eigenvalues can be bounded by threshold graphs leads
to Brouwer's conjecture. It is related to (and motivated by)
the Grone-Merris conjecture  \cite{Haemers20102214} and states that for any graph
$G$ on $n$ vertices and $m$ edges,
\begin{equation}\label{Brouwer}
\sum_{i=1}^k\lambda_i(G) \leq m + {k+1 \choose 2}.
\end{equation}
One may ask whether the bound given by the Grone-Merris-Bai theorem is
sharper than Brouwer's conjecture, because it uses more detailed
information from the graph. Indeed, it has been shown in \cite{Mayank} that for split
graphs this is the case. However, more generally, it is shown that
there is a $k$ such that the $k$-th inequality of Brouwer's conjecture
is sharper than the $k$-th Grone-Merris inequality if and only if the
graph is non-split. Brouwer's conjecture remains unproven to this
date.

It turns out (see \cref{S:brouwer})
that Brouwer's conjecture is, in fact, equivalent to spectral
threshold dominance.
\begin{Theorem}\label{Th:brouwer}
  A graph $G$ satisfies Brouwer's conjecture if and only if it is
  spectrally threshold dominated.
\end{Theorem}
Quite likely this relation to threshold graphs has been part of the motivation
for Brouwer's conjecture. Recognizing this equivalence also opens
the door to previous, rather different proofs of \cref{Th:split},
\cref{Th:preserve}, and \cref{Th:cographs} by \cite{Mayank}, who
proved that Brouwer's conjecture holds in these cases.

Establishing Brouwer's conjecture would prove the Laplacian energy
conjecture in the non connected case. The requirement of spectral
threshold dominance is, however, stronger than needed for the
Laplacian energy conjecture. A counterexample for Brouwer's
conjecture might not be sufficient to disprove the Laplacian energy
conjecture. On the other hand, a counterexample for the non connected
Laplacian energy conjecture would immediately disprove the spectral
threshold dominance conjecture and thus also Brouwer's conjecture.

\section{Spectral threshold dominance and Laplacian energy}\label{S:specdom}

Before embarking on the proof of \cref{Th:specdom} we illustrate
spectral threshold dominance by an example involving cycles (these are
neither cographs nor split graphs). Note that in the definition of
spectral threshold dominance, the threshold graph is allowed to depend
on $G$ and $k$. In the example, we are able to supply a single
threshold graph that spectrally dominates $G$ for all values of
$k$.

Our constructions here and later are inspired by the
characterization of threshold graphs by their Ferrers (or Young)
diagram (see for example \cite{MerrisRoby2005}) of their nonincreasing
degree sequence. Given a degree sequence $d_1 \geq d_2 \geq \cdots \geq d_n$ of
a graph $G$, its Ferrers diagram consists of $n$ rows, where row $i$
displays $d_i$ boxes aligned on the left. The conjugate degree sequence $d_i^*$
may be read off from the Ferrers diagram by simply counting the number of boxes
in the $i$-the column. The diagonal width of the degree sequence $f=\max\{i\in N\colon
d_i\ge i\}$ is called the \emph{trace} of the diagram.

For threshold graphs, as $d_i^* = d_i+1$, it means that in each row $i \in [f]$ ,
the number of boxes on and to the right of the diagonal box is the same as the number
of boxes in column $i$ below the diagonal. More generally, the shape described by
the boxes on and above the diagonal is exactly the transpose of the shape of the
boxes below the diagonal. The $f$ boxes on the diagonal will be displayed in
black.

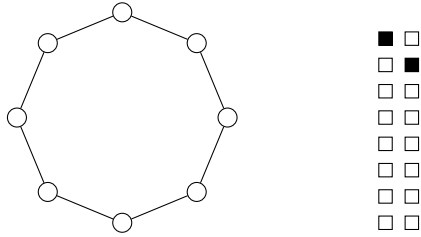
\begin{figure}[h!]\begin{center}
\begin{tikzpicture}
  [fill=black!100,scale=0.7,auto=left,every node/.style={circle,scale=0.7}]
  \foreach \i in {1,...,8}{
    \node[draw,circle] (\i) at ({360/8 * (\i - 1)+90}:2) {};} % :2 is the radius
\draw (1) -- (2) -- (3) -- (4) -- (5) -- (6) -- (7) -- (8) -- (1);

\node at ( 5,1.5) [rectangle,draw,fill] {};
\node at ( 5,.5) [rectangle,draw] {};
\node at ( 5,1) [rectangle,draw] {};
\node at ( 5,0) [rectangle,draw] {};
\node at ( 5,-.5) [rectangle,draw] {};
\node at ( 5,-1) [rectangle,draw] {};
\node at ( 5,-1.5) [rectangle,draw] {};
\node at ( 5,-2) [rectangle,draw] {};

\node at ( 5.5,1.5) [rectangle,draw] {};
\node at ( 5.5,.5) [rectangle,draw] {};
\node at ( 5.5,1) [rectangle,draw,fill] {};
\node at ( 5.5,0) [rectangle,draw] {};
\node at ( 5.5,-.5) [rectangle,draw] {};
\node at ( 5.5,-1) [rectangle,draw] {};
\node at ( 5.5,-1.5) [rectangle,draw] {};
\node at ( 5.5,-2) [rectangle,draw] {};
\end{tikzpicture}
    \caption{ $C_8$ and its Ferrers diagram}
    \label{fig1ex}
\end{center}\end{figure}

\begin{Example} Figure \ref{fig1ex} depicts the cycle $C_8$ on 8
  vertices and its Ferrers diagram, while Figure~\ref{fig2ex} shows
  a threshold graph that
  spectrally dominates $C_8$ for all $k \in \{1,\ldots,8\}$. Indeed
  the spectrum of $C_8$ is the multiset  $\{ 4, 2+\sqrt{2},2+\sqrt{2},
  2, 2, 2-\sqrt{2},2-\sqrt{2}, 0\}$ whose partial sums are $4, 6 +
  \sqrt{2}, 8 + 2\sqrt{2}, 10 + 2\sqrt{2}, 12 + 2\sqrt{2}, 14+
  \sqrt{2}, 16,16$, whereas the partial sums  for the eigenvalues of
  the threshold graph of Figure \ref{fig2ex} are
  $5,10,14,16,16,16,16,16$. Notice that the resulting threshold graph is disconnected.

\begin{figure}[h!]\begin{center}
\begin{tikzpicture}
[fill=black!150,scale=0.7,auto=left,every node/.style={circle,scale=0.7}]
  \foreach \i in {1,...,8}{
    \node[draw,circle] (\i) at ({360/8 * (\i - 1)+90}:2) {};} % :2 is the radius
\draw (1) -- (2) -- (3) -- (4);
\draw (2) -- (5);
\draw (2) -- (4);
\draw (1) -- (3);
\draw (1) -- (4);
\draw (1) -- (5);

\node at ( 4, 1.5) [rectangle,draw,fill] {};
\node at ( 4.5, 1.5) [rectangle,draw] {};
\node at ( 5, 1.5) [rectangle,draw] {};
\node at ( 5.5, 1.5) [rectangle,draw] {};

\node at ( 4, 1) [rectangle,draw] {};
\node at ( 4.5, 1) [rectangle,draw,fill] {};
\node at ( 5, 1) [rectangle,draw] {};
\node at ( 5.5, 1) [rectangle,draw] {};

\node at ( 5, 0.5) [rectangle,draw,fill] {};
\node at ( 4, 0.5) [rectangle,draw] {};
\node at ( 4.5,0.5) [rectangle,draw] {};

\node at ( 5, 0) [rectangle,draw] {};
\node at ( 4, 0) [rectangle,draw] {};
\node at ( 4.5, 0) [rectangle,draw] {};

\node at ( 4, -.5) [rectangle,draw] {};
\node at ( 4.5, -.5) [rectangle,draw] {};
\end{tikzpicture}
 \caption{A spectrally threshold dominant graph of $C_8$}
     \label{fig2ex}
\end{center}\end{figure}
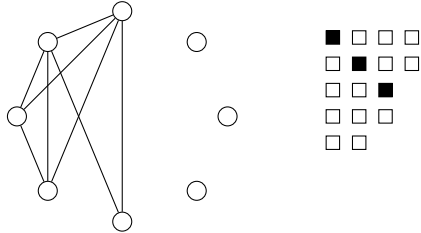

This procedure can be generalized for a general cycle $C_n$ on $n$
vertices. Indeed consider $C_n$ with $n\geq 8$ vertices (hence with
$m=n \geq 8$ edges).  Let $h=\lfloor \sqrt{2n} \rfloor$. We observe
that $ 2n- (h^2-h) \geq h$ and $2n-(h^2+h) \leq h+1$. Define $T$ to be  the
threshold graph whose Ferrers diagram has Durfee square (trace)
$f=h-1$, if $2n-(h^2+h) < 0$, otherwise  $T$ will have trace $f=h$. The
remaining boxes will be placed as a last $(f+1)$th-column (half of
them) and the corresponding last $(f+2)$th-row in the Ferrers diagram.  As the
Laplacian spectrum of $C_n$ is the set  $\{ 2- 2\cos
\frac{2\pi}{n}i\colon i=1, \ldots, n\}$ we observe that for
$i=1,\ldots, f$, $\lambda_i(C_n) \leq 4 = \sqrt{2\cdot 8} = \lfloor
\sqrt{2\cdot 8} \rfloor \leq \lfloor \sqrt{2\cdot n} \rfloor \leq
\lambda_i(T)$. Hence it holds that $\sum_{i=1}^k \lambda_i(C_n) \leq
4\cdot k \leq h\cdot k \leq \sum_{i=1}^k \lambda_i(T)$, for $k=1, \ldots, f$. If $k\geq f$, we observe that $\sum_{i=1}^k \lambda_i(T) = 2n =2m \geq \sum_{i=1}^k \lambda_i(C_n)$.
\end{Example}

The Laplacian energy of a graph $G$ is actually  fully determined by the sum
of the $k$ eigenvalues whose values exceed the average degree
$\tfrac{2m}n$. For providing a threshold graph $T$ on the same number of
nodes and edges with the same or higher Laplacian energy it suffices
to find one with $\sum_{i=1}^k\lambda_i(T)\ge\sum_{i=1}^k\lambda_i(G)$
for this specific $k$, as proved in the following lemma.

\begin{Lemma}\label{L:dominated}
Let $G$ be a graph on $n$ nodes with $m$ edges and conjugate degree sequence $d_i^*$, $i=1,\dots,n$ and let $k\in\{1,\dots,n\}$ be the index satisfying
$\lambda_k(G)> \tfrac{2m}{n}\ge\lambda_{k+1}(G)$. Any threshold graph $T$ on $n$ nodes with $m$ edges satisfying $\sum_{i=1}^kd_i^*(T)\ge \sum_{i=1}^k\lambda_i(G)$ also satisfies $LE(T)\ge LE(G)$.
\end{Lemma}
\begin{Proof}
By $\sum_{i=0}^n\lambda_i(G)=2m$ there holds $\sum_{i=k+1}^n\lambda_i(G)=2m-\sum_{i=1}^k\lambda_i(G)$ and therefore
\begin{eqnarray*}
LE(G)&=&\sum_{i=1}^n|\lambda_i(G)-\tfrac{2m}n|\\
&=&\sum_{i=1}^k
(\lambda_i(G)-\tfrac{2m}n)+\sum_{i=k+1}^n(\tfrac{2m}{n}-\lambda_i(G))\\
&=&2\sum_{i=1}^k \lambda_i(G)-2m-k\tfrac{2m}{n}+(n-k)\tfrac{2m}{n}\\
&=&2\sum_{i=1}^k(\lambda_i(G)-\tfrac{2m}n)\\
&\le &2\sum_{i=1}^k(d_i^*(T)-\tfrac{2m}n)\\
&=&2\sum_{i=1}^k(\lambda_i(T)-\tfrac{2m}n)\\
&\le&2\sum_{i\in\{j\colon \lambda_j(T)>\tfrac{2m}{n}\}}(\lambda_i(T)-\tfrac{2m}n)=LE(T).
\end{eqnarray*}
The last equation follows by repeating the initial arguments on $G$ for $T$.
\end{Proof}
If a graph is spectrally threshold dominated, an appropriate threshold
graph $T_k$ exists for all $k$, in particular for the $k$ required in
\cref{L:dominated}. This proves \cref{Th:specdom}.

\section{Split graphs are spectrally threshold dominated}\label{S:split}

Recall that $G$ is a split graph if its set of vertices can be
partitioned in two sets $A$ and $B$ such that $A$ induces a clique in
$G$ and $B$ does not contain any edge. The key for proving that split graphs are
spectrally threshold dominated is the characterization of split graphs
and threshold graphs by their Ferrers diagram. Split graphs have the
same number of boxes above and on the diagonal as below the diagonal.
Threshold graphs are special split graphs in that the shape below is
the transposed of the shape above and on the diagonal.
This forms the basis of the proof of the following lemma, which
directly establishes \cref{Th:split}.

\begin{Lemma} Given a split graph $G$ on $n$ nodes with $m$ edges and $k\in\{1,\dots,n\}$, there is a threshold graph $T$ on $n$ nodes with $m$ edges so that $\sum_{i=1}^k\lambda_i(T)\ge\sum_{i=1}^k\lambda_i(G).$
\end{Lemma}
\begin{Proof}
Let $f(G)=\max\{i\colon d_i(G)\ge i\}$ be the trace of the Ferrers diagram for $G$. We discern the cases $k<f(G)$ and $k\ge f(G)$.

\begin{figure}[h!]\begin{center}
\begin{tikzpicture}
[fill=black!150,scale=0.7,auto=left,every node/.style={circle,scale=0.7}]

\node at ( 0, 1.5) [rectangle,draw,fill] {};
\node at ( .5, 1.5) [rectangle,draw] {};
\node at ( 1, 1.5) [rectangle,draw] {};
\node at ( 1.5, 1.5) [rectangle,draw] {};

\node at ( 0, 1) [rectangle,draw] {};
\node at ( .5, 1) [rectangle,draw,fill] {};
\node at ( 1, 1) [rectangle,draw] {};
\node at ( 1.5, 1) [rectangle,draw] {};

\node at ( 0, .5) [rectangle,draw] {};
\node at ( .5, .5) [rectangle,draw] {};
\node at ( 1, .5) [rectangle,draw,fill=red!150,pattern=north east lines] {};

\node at ( 0, 0) [rectangle,draw] {};
\node at ( .5, 0) [rectangle,draw] {};

\node at ( 0, -.5) [rectangle,draw] {};
\node at ( .5, -.5) [rectangle,draw] {};

\node at ( 0, -1) [rectangle,draw,] {};

\node at ( 4, 1.5) [rectangle,draw,fill] {};
\node at ( 4.5, 1.5) [rectangle,draw] {};
\node at ( 5, 1.5) [rectangle,draw] {};
\node at ( 5.5, 1.5) [rectangle,draw] {};
\node at ( 6, 1.5) [rectangle,draw,fill=red!150,pattern=north east lines] {};

\node at ( 4, 1) [rectangle,draw] {};
\node at ( 4.5, 1) [rectangle,draw,fill] {};
\node at ( 5, 1) [rectangle,draw] {};
\node at ( 5.5, 1) [rectangle,draw] {};

\node at ( 4, .5) [rectangle,draw] {};
\node at ( 4.5, .5) [rectangle,draw] {};

\node at ( 4, 0) [rectangle,draw] {};
\node at ( 4.5, 0) [rectangle,draw] {};

\node at ( 4, -.5) [rectangle,draw] {};
\node at ( 4.5, -.5) [rectangle,draw] {};

\node at ( 4, -1) [rectangle,draw] {};

\end{tikzpicture}
 \caption{A transform for $k < f(G)$.}
     \label{smallk}
\end{center}\end{figure}
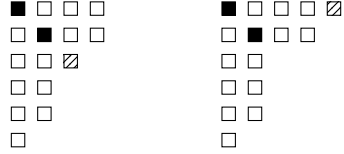

\underline{$k<f(G)$:} leaving the boxes below the diagonal unchanged and by copying its  shape  in transposed form to the part above the diagonal we obtain a diagram uniquely defining a threshold graph $T$ with the property $\sum_{i=1}^k\lambda_i(T) = \sum_{i=1}^kd_i^*(T)= \sum_{i=1}^kd_i^*(G)\ge\sum_{i=1}^k\lambda_i(G)$, where the last inequality follows from the Grone-Merris-Bai theorem. We refer to Figure \ref{smallk} for an illustration of the graph transform. In the example $k=2$ and the hatched box is moved.

\begin{figure}[h!]\begin{center}
\begin{tikzpicture}
[fill=black!150,scale=0.7,auto=left,every node/.style={circle,scale=0.7}]

\node at ( 0, 1.5) [rectangle,draw,fill] {};
\node at ( .5, 1.5) [rectangle,draw] {};
\node at ( 1, 1.5) [rectangle,draw] {};
\node at ( 1.5, 1.5) [rectangle,draw] {};
\node at ( 2, 1.5) [rectangle,draw,fill=blue!150,pattern=north east lines] {};
\node at ( 2.5, 1.5) [rectangle,draw,fill=blue!150,pattern=north east lines] {};

\node at ( 0, 1) [rectangle,draw] {};
\node at ( .5, 1) [rectangle,draw,fill] {};
\node at ( 1, 1) [rectangle,draw] {};

\node at ( 0, .5) [rectangle,draw] {};
\node at ( .5, .5) [rectangle,draw] {};
\node at ( 1, .5) [rectangle,draw,fill] {};

\node at ( 0, 0) [rectangle,draw] {};
\node at ( .5, 0) [rectangle,draw] {};
\node at ( 1, 0) [rectangle,draw] {};

\node at ( 0, -.5) [rectangle,draw] {};
\node at ( .5, -.5) [rectangle,draw] {};

\node at ( 0, -1) [rectangle,draw,fill=red!150,pattern=north west lines] {};

\node at ( 4, 1.5) [rectangle,draw,fill] {};
\node at ( 4.5, 1.5) [rectangle,draw] {};
\node at ( 5, 1.5) [rectangle,draw] {};
\node at ( 5.5, 1.5) [rectangle,draw] {};

\node at ( 4, 1) [rectangle,draw] {};
\node at ( 4.5, 1) [rectangle,draw,fill] {};
\node at ( 5, 1) [rectangle,draw] {};
\node at ( 5.5, 1) [rectangle,draw,fill=blue!150,pattern=north east lines] {};

\node at ( 4, .5) [rectangle,draw] {};
\node at ( 4.5, .5) [rectangle,draw] {};
\node at ( 5, .5) [rectangle,draw,fill] {};
\node at ( 5.5, .5) [rectangle,draw,fill=blue!150,pattern=north east lines] {};

\node at ( 4, 0) [rectangle,draw] {};
\node at ( 4.5, 0) [rectangle,draw] {};
\node at ( 5, 0) [rectangle,draw] {};

\node at ( 4, -.5) [rectangle,draw] {};
\node at ( 4.5, -.5) [rectangle,draw] {};
\node at ( 5, -.5) [rectangle,draw,fill=red!150,pattern=north west lines] {};

\end{tikzpicture}
 \caption{A transform for $k \ge f(G)$.}
     \label{largek}
\end{center}\end{figure}

\underline{$k\ge f(G)$:} construct the desired threshold graph $T$ by
filling up the diagram above and on the diagonal in columnwise order
by the $m$ boxes, but only up to and including row $f(G)$ (and in
rowwise order below the diagonal up to column $f(G)$ for the
transposed shape). Thus $f(T)=f(G)$. Figure \ref{largek} shows the
threshold graph $T$ for a particular graph $G$. In the example the
hatched boxes are moved. The construction never moves boxes across the diagonal.

First consider the case $k=f(G)$. Because the number of boxes below the diagonal
is the same for $G$ and $T$, we obtain
\begin{equation}
  \label{eq:same_sum}
\sum_{i=1}^{f(G)}\lambda_i(T) =\sum_{i=1}^{f(G)}d_i^*(T) =m+f(G)(f(G)+1)/2 =
\sum_{i=1}^{f(G)}d_i^*(G).
\end{equation}
The claim follows from
$\sum_{i=1}^kd_i^*(G) \ge \sum_{i=1}^k\lambda_i(G)$ by the
Grone-Merris-Bai theorem.

Finally, for $k>f(G)$ observe that for $j \in \{f(G)+1, \ldots, n\}$
there holds
$\sum_{i=f(G)+1}^j d_i^*(T) \geq\sum_{i=f(G)+1}^j  d_i^*(G)$, because
in $T$ the boxes have been rearranged to maximally fill up the first
columns after column $f(G)$. Thus, equation \eqref{eq:same_sum} and the
Grone-Merris-Bai theorem yield $\sum_{i=1}^kd_i^*(T) \ge \sum_{i=1}^kd_i^*(G) \ge \sum_{i=1}^k\lambda_i(G)$. \end{Proof}

Note that the construction of the proof for $k\ge f(G)$ may generate a
threshold graph that is not connected even if $G$ is. Indeed, at this
point we do not know how to construct for a general connected split
graph $G$ and given $k$ a spectrally dominating connected threshold
graph $T_k$.

\section{Disjoint unions and complements preserve spectral threshold dominance}
\label{S:preserve}

In order to increase the class of spectrally threshold dominated graphs a bit further, we consider taking the union and complements of spectrally threshold dominated graphs.

\begin{Lemma}\label{L:union}
Let $G$ be a (disjoint) union of spectrally threshold dominated graphs with $n$ nodes and $m$ edges and let $k\in\{1,\dots,n\}$. There is a threshold graph $T$ on $n$ nodes and $m$ edges so that $\sum_{i=1}^k\lambda_i(T)\ge\sum_{i=1}^k\lambda_i(G).$
\end{Lemma}

Before proving that the disjoint union of spectrally threshold dominated
graphs $G_j=(N_j,\bar E_j)$, $j=1,\dots,h$ is spectrally threshold dominated
we first explain the main idea of the approach. For each $j=1,\dots,h$ put
$n_j=|N_j|$ and suppose the eigenvalues of $G_j$ are sorted non increasingly
$\lambda_1(G_j)\ge \dots\ge\lambda_{n_j}(G_j)$.  The spectrum of the union
$G=\bigcup_{j=1}^hG_j$ of the graphs is the union of the spectra of the
graphs. Let the $k$ largest eigenvalues of $G$ arise from taking the $k_j$
largest eigenvalues of $G_j$ for $j=1,\dots,h$, then $k=\sum_{j=1}^h k_j$. If
some graphs have identical eigenvalues, it does not matter which of them are
selected. Because each $G_j$ is spectrally dominated, there is a threshold
graph $T_j=(N_j,E_j)$ with $|E_j|=|\bar E_j|$ so that its non decreasingly
sorted dual degrees $d_{(j,i)}^*:=d_i^*(T_j)$, which equal the eigenvalues
$\lambda_i(T_j)=d_{(j,i)}^*$ satisfy $\sum_{i=1}^{k_j}d_{(j,i)}^*\ge
\sum_{i=1}^{k_j}\lambda_i(G_j)$. Therefore it suffices to consider
$G_j=T_j$, take the $k$ largest with respect to the sorting of the
$d_{(j,i)}^*=\lambda_i(T_j)$ (this only increases the sum) and construct
a threshold graph $T$ with the same number of nodes and edges as
$G=\bigcup_{j=1}^hT_j$ so that its non increasing dual degrees $d_i^*(T)$
satisfy $\sum_{i=1}^{k}d^*_i(T) \ge \sum_{j=1}^h\sum_{i=1}^{k_j}d_{(j,i)}^*$. This
$T$ arises by the following simple construction of its Ferrers diagram. In
column $i$ the boxes below the diagonal are obtained by concatenating over all
$j=1,\dots,h$ the boxes of column $i$ that are below the diagonal of the
diagram of $T_j$. The part above the diagonal is the corresponding transposed
structure. In order to see that this graph actually dominates the sum of the
$k$ largest eigenvalues of $G$, an intermediate construction step is needed,
where the columns of all the diagrams are arranged side by side and sorted by
their height. Then, as illustrated in the example below, $T$ arises from
this diagram by only shifting boxes from higher column indices to lower column
indices, so the sum of the boxes in the first $k$ columns has to increase.

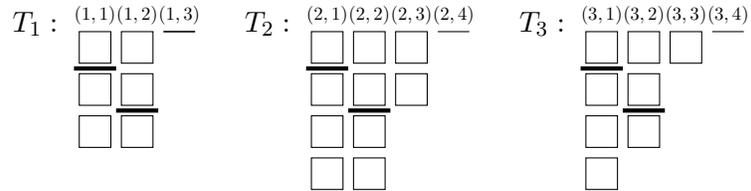
\begin{center}
$T_1:$ \begin{tikzpicture}[baseline=0ex,scale=0.7]
\foreach \name in {1}
{
\foreach \x in {1}
{
\draw (.8*\x-.8+.3,0) node[above,scale=.7]{\small$(\name,\x)$};
\draw[ultra thick] (.8*\x-.9,-.8*\x+.1)--(.8*\x-.8+.7,-.8*\x+.1);
\foreach \y in {0,1,2}
{
 \draw (.8*\x-.8,-.8*\y) rectangle (.8*\x-.8+ .6,-.8 *\y -.6);
}
}
\foreach \x in {2}
{
\draw (.8*\x-.8+.3,0) node[above,scale=.7]{\small$(\name,\x)$};
\draw[ultra thick] (.8*\x-.9,-.8*\x+.1)--(.8*\x-.8+.7,-.8*\x+.1);
\foreach \y in {0,1,2}
{
 \draw (.8*\x-.8,-.8*\y) rectangle (.8*\x-.8 + .6,-.8 *\y -.6);
}
}
\foreach \x in {3}
{
\draw (.8*\x-.8+.3,0) node[above,scale=.7]{\small$(\name,\x)$};
\foreach \y in {0,1}
\foreach \y in {0}
{
 \draw (.8*\x-.8,-.8*\y) -- (.8*\x-.8 + .6,-.8 *\y);
}
}
}
\end{tikzpicture}
\quad $T_2:$  \begin{tikzpicture}[baseline=0ex,scale=.7]
\foreach \name in {2}
{
\foreach \x in {1}
{
\draw (.8*\x-.8+.3,0) node[above,scale=.7]{\small$(\name,\x)$};
\draw[ultra thick] (.8*\x-.9,-.8*\x+.1)--(.8*\x-.8+.7,-.8*\x+.1);
\foreach \y in {0,1,2,3}
{
 \draw (.8*\x-.8,-.8*\y) rectangle (.8*\x-.8+ .6,-.8 *\y -.6);
}
}
\foreach \x in {2}
{
\draw (.8*\x-.8+.3,0) node[above,scale=.7]{\small$(\name,\x)$};
\draw[ultra thick] (.8*\x-.9,-.8*\x+.1)--(.8*\x-.8+.7,-.8*\x+.1);
\foreach \y in {0,1,2,3}
{
 \draw (.8*\x-.8,-.8*\y) rectangle (.8*\x-.8 + .6,-.8 *\y -.6);
}
}
\foreach \x in {3}
{
\draw (.8*\x-.8+.3,0) node[above,scale=.7]{\small$(\name,\x)$};
\foreach \y in {0,1}
{
 \draw (.8*\x-.8,-.8*\y) rectangle (.8*\x-.8 + .6,-.8 *\y -.6);
}
}
\foreach \x in {4}
{
\draw (.8*\x-.8+.3,0) node[above,scale=.7]{\small$(\name,\x)$};
\foreach \y in {0}
{
 \draw (.8*\x-.8,-.8*\y) -- (.8*\x-.8 + .6,-.8 *\y);
}
}
}
\end{tikzpicture}
\quad $T_3:$
\begin{tikzpicture}[baseline=0ex,scale=.7]
\foreach \name in {3}
{
\foreach \x in {1}
{
\draw (.8*\x-.8+.3,0) node[above,scale=.7]{\small$(\name,\x)$};
\draw[ultra thick] (.8*\x-.9,-.8*\x+.1)--(.8*\x-.8+.7,-.8*\x+.1);
\foreach \y in {0,1,2,3}
{
 \draw (.8*\x-.8,-.8*\y) rectangle (.8*\x-.8+ .6,-.8 *\y -.6);
}
}
\foreach \x in {2}
{
\draw (.8*\x-.8+.3,0) node[above,scale=.7]{\small$(\name,\x)$};
\draw[ultra thick] (.8*\x-.9,-.8*\x+.1)--(.8*\x-.8+.7,-.8*\x+.1);
\foreach \y in {0,1,2}
{
 \draw (.8*\x-.8,-.8*\y) rectangle (.8*\x-.8 + .6,-.8 *\y -.6);
}
}
\foreach \x in {3}
{
\draw (.8*\x-.8+.3,0) node[above,scale=.7]{\small$(\name,\x)$};
\foreach \y in {0}
{
 \draw (.8*\x-.8,-.8*\y) rectangle (.8*\x-.8 + .6,-.8 *\y -.6);
}
}
\foreach \x in {4}
{
\draw (.8*\x-.8+.3,0) node[above,scale=.7]{\small$(\name,\x)$};
\foreach \y in {0}
{
 \draw (.8*\x-.8,-.8*\y) -- (.8*\x-.8 + .6,-.8 *\y);
}
}
}
\end{tikzpicture}\\
\end{center}

\begin{center}
sorted eigenvalues: \begin{tikzpicture}[baseline=0ex,scale=.7]
\foreach \x in {1}
{
\draw (.8*\x-.8+.3,0) node[above,scale=.7]{\small$(2,1)$};
\draw[ultra thick] (.8*\x-.9,-.8*1+.1)--(.8*\x-.8+.7,-.8*1+.1);
\foreach \y in {0,1,2,3}
{
 \draw (.8*\x-.8,-.8*\y) rectangle (.8*\x-.8+ .6,-.8 *\y -.6);
}
\draw[dotted,very thick] (.8*\x-.8-.03,-.8*1+.03) rectangle (.8*\x-.8+.6+.03,-.8*3-.6-.03);
\draw[->] (.8*\x-.8+.3,-.8*4)--(.8*\x-.8+.3,-.8*4-.3) node[below]{1};
}
\foreach \x in {2}
{
\draw (.8*\x-.8+.3,0) node[above,scale=.7]{\small$(3,1)$};
\draw[ultra thick] (.8*\x-.9,-.8*1+.1)--(.8*\x-.8+.7,-.8*1+.1);
\foreach \y in {0,1,2,3}
{
 \draw (.8*\x-.8,-.8*\y) rectangle (.8*\x-.8 + .6,-.8 *\y -.6);
}
\draw[dotted,very thick] (.8*\x-.8-.03,-.8*1+.03) rectangle (.8*\x-.8+.6+.03,-.8*3-.6-.03);
\draw[->] (.8*\x-.8+.3,-.8*4)--(.8*\x-.8+.3,-.8*4-.3) node[below]{1};
}
\foreach \x in {3}
{
\draw (.8*\x-.8+.3,0) node[above,scale=.7]{\small$(2,2)$};
\draw[ultra thick] (.8*\x-.9,-.8*2+.1)--(.8*\x-.8+.7,-.8*2+.1);
\foreach \y in {0,1,2,3}
{
 \draw (.8*\x-.8,-.8*\y) rectangle (.8*\x-.8 + .6,-.8 *\y -.6);
}
\draw[dotted,very thick] (.8*\x-.8-.03,-.8*2+.03) rectangle (.8*\x-.8+.6+.03,-.8*3-.6-.03);
\draw[->] (.8*\x-.8+.3,-.8*4)--(.8*\x-.8+.3,-.8*4-.3) node[below]{2};
\draw (.8*\x-.8+.3,-.8*1-.3) node{$\leftarrow$};
}
\foreach \x in {4}
{
\draw (.8*\x-.8+.3,0) node[above,scale=.7]{\small$(1,1)$};
\draw[ultra thick] (.8*\x-.9,-.8*1+.1)--(.8*\x-.8+.7,-.8*1+.1);
\foreach \y in {0,1,2}
{
 \draw (.8*\x-.8,-.8*\y) rectangle (.8*\x-.8 + .6,-.8 *\y -.6);
}
\draw[dotted,very thick] (.8*\x-.8-.03,-.8*1+.03) rectangle (.8*\x-.8+.6+.03,-.8*2-.6-.03);
\draw[->] (.8*\x-.8+.3,-.8*3)--(.8*\x-.8+.3,-.8*3-.3) node[below]{1};
}
\foreach \x in {5}
{
\draw (.8*\x-.8+.3,0) node[above,scale=.7]{\small$(1,2)$};
\draw[ultra thick] (.8*\x-.9,-.8*2+.1)--(.8*\x-.8+.7,-.8*2+.1);
\foreach \y in {0,1,2}
{
 \draw (.8*\x-.8,-.8*\y) rectangle (.8*\x-.8 + .6,-.8 *\y -.6);
}
\draw[dotted,very thick] (.8*\x-.8-.03,-.8*2+.03) rectangle (.8*\x-.8+.6+.03,-.8*2-.6-.03);
\draw[->] (.8*\x-.8+.3,-.8*3)--(.8*\x-.8+.3,-.8*3-.3) node[below]{2};
\draw (.8*\x-.8+.3,-.8*1-.3) node{$\leftarrow$};
}
\foreach \x in {6}
{
\draw (.8*\x-.8+.3,0) node[above,scale=.7]{\small$(3,2)$};
\draw[ultra thick] (.8*\x-.9,-.8*2+.1)--(.8*\x-.8+.7,-.8*2+.1);
\foreach \y in {0,1,2}
{
 \draw (.8*\x-.8,-.8*\y) rectangle (.8*\x-.8 + .6,-.8 *\y -.6);
}
\draw[dotted,very thick] (.8*\x-.8-.03,-.8*2+.03) rectangle (.8*\x-.8+.6+.03,-.8*2-.6-.03);
\draw[->] (.8*\x-.8+.3,-.8*3)--(.8*\x-.8+.3,-.8*3-.3) node[below]{2};
\draw (.8*\x-.8+.3,-.8*1-.3) node{$\leftarrow$};
}
\foreach \x in {7}
{
\draw (.8*\x-.8+.3,0) node[above,scale=.7]{\small$(2,2)$};
\draw[ultra thick] (.8*\x-.9,-.8*2+.1)--(.8*\x-.8+.7,-.8*2+.1);
\foreach \y in {0,1}
{
 \draw (.8*\x-.8,-.8*\y) rectangle (.8*\x-.8 + .6,-.8 *\y -.6);
}
\draw (.8*\x-.8+.3,-.8*1-.3) node{$\leftarrow$};
}
\foreach \x in {8}
{
\draw (.8*\x-.8+.3,0) node[above,scale=.7]{\small$(3,3)$};
\foreach \y in {0}
{
 \draw (.8*\x-.8,-.8*\y) rectangle (.8*\x-.8 + .6,-.8 *\y -.6);
}
}
\foreach \x in {9}
{
\draw (.8*\x-.8+.3,0) node[above,scale=.7]{\small$(1,3)$};
\foreach \y in {0}
{
 \draw (.8*\x-.8,-.8*\y) -- (.8*\x-.8 + .6,-.8 *\y);
}
}
\foreach \x in {10}
{
\draw (.8*\x-.8+.3,0) node[above,scale=.7] {\small$(2,4)$};
\foreach \y in {0}
{
 \draw (.8*\x-.8,-.8*\y) -- (.8*\x-.8 + .6,-.8 *\y);
}
}
\foreach \x in {11}
{
\draw (.8*\x-.8+.3,0) node[above,scale=.7]{\small$(3,4)$};
\foreach \y in {0}
{
 \draw (.8*\x-.8,-.8*\y) -- (.8*\x-.8 + .6,-.8 *\y);
}
}
\end{tikzpicture}\\
\end{center}

\begin{center}
final $T$:
\begin{tikzpicture}[baseline=0ex,scale=.7]
\foreach \x in {1}
{
\draw[ultra thick] (.8*\x-.9,-.8*\x+.1)--(.8*\x-.8+.7,-.8*\x+.1);
\foreach \y in {0,...,8}
{
 \draw (.8*\x-.8,-.8*\y) rectangle (.8*\x-.8+ .6,-.8 *\y -.6);
}
\draw[dotted,very thick] (.8*\x-.8-.03,-.8*1+.03) rectangle (.8*\x-.8+.6+.03,-.8*3-.6-.03);
\draw[dotted,very thick] (.8*\x-.8-.03,-.8*4+.03) rectangle (.8*\x-.8+.6+.03,-.8*6-.6-.03);
\draw[dotted,very thick] (.8*\x-.8-.03,-.8*7+.03) rectangle (.8*\x-.8+.6+.03,-.8*8-.6-.03);
}
\foreach \x in {2}
{
\draw[ultra thick] (.8*\x-.9,-.8*\x+.1)--(.8*\x-.8+.7,-.8*\x+.1);
\foreach \y in {0,...,5}
{
 \draw (.8*\x-.8,-.8*\y) rectangle (.8*\x-.8 + .6,-.8 *\y -.6);
}
\draw[dotted,very thick] (.8*\x-.8-.03,-.8*2+.03) rectangle (.8*\x-.8+.6+.03,-.8*3-.6-.03);
\draw[dotted,very thick] (.8*\x-.8-.03,-.8*4+.03) rectangle (.8*\x-.8+.6+.03,-.8*4-.6-.03);
\draw[dotted,very thick] (.8*\x-.8-.03,-.8*5+.03) rectangle (.8*\x-.8+.6+.03,-.8*5-.6-.03);
}
\foreach \x in {3}
{
\foreach \y in {0,1}
{
 \draw (.8*\x-.8,-.8*\y) rectangle (.8*\x-.8 + .6,-.8 *\y -.6);
}
}
\foreach \x in {4}
{
\foreach \y in {0,1}
{
 \draw (.8*\x-.8,-.8*\y) rectangle (.8*\x-.8 + .6,-.8 *\y -.6);
}
}
\foreach \x in {5}
{
\foreach \y in {0,1}
{
 \draw (.8*\x-.8,-.8*\y) rectangle (.8*\x-.8 + .6,-.8 *\y -.6);
}
}
\foreach \x in {6}
{
\foreach \y in {0}
{
 \draw (.8*\x-.8,-.8*\y) rectangle (.8*\x-.8 + .6,-.8 *\y -.6);
}
}
\foreach \x in {7}
{
\foreach \y in {0}
{
 \draw (.8*\x-.8,-.8*\y) rectangle (.8*\x-.8 + .6,-.8 *\y -.6);
}
}
\foreach \x in {8}
{
\foreach \y in {0}
{
 \draw (.8*\x-.8,-.8*\y) rectangle (.8*\x-.8 + .6,-.8 *\y -.6);
}
}
\foreach \x in {9}
{
\foreach \y in {0}
{
 \draw (.8*\x-.8,-.8*\y) -- (.8*\x-.8 + .6,-.8 *\y);
}
}
\foreach \x in {10}
{
\foreach \y in {0}
{
 \draw (.8*\x-.8,-.8*\y) -- (.8*\x-.8 + .6,-.8 *\y);
}
}
\foreach \x in {11}
{
\foreach \y in {0}
{
 \draw (.8*\x-.8,-.8*\y) -- (.8*\x-.8 + .6,-.8 *\y);
}
}
\end{tikzpicture}
\end{center}
Unfortunately, the algebraic formulation of this geometrically intuitive
observation requires an algorithmic description of the intermediate step
and is therefore notationally  more involved.

\begin{Proof}(Lemma \ref{L:union})
We assume that $G=\bigcup_{j=1}^{h}G_j$ with each $G_j$ spectrally
threshold dominated. Let the first $k$ eigenvalues of $G$ consist of
the first $k_j$ eigenvalues of $G_j$, $j=1,\dots,h$ with
$\sum_{j=1}^hk_j=k$. For each $G_j$ there is threshold graph $T_j$
so that $\sum_{i=1}^{k_j} \lambda_i(T_j) \ge
\sum_{i=1}^{k_j}\lambda_i(G_j).$ As seen above, it suffices to
prove the  result under the assumption that each $G_j$ is a threshold
graph $T_j=(N_j,E_j)$ with $n_j=|N_j|$, $m_j=|E_j|$ so that
$\sum_{j=1}^hn_j=n$ and $\sum_{j=1}^hm_j=m$.   Put $H=\{(j,i): j\in\{1,\dots,h\},i\in\{1,\dots,n_j\}\}$ with index
  $(j,i)$ representing $d_{(j,i)}^*=d_i^*(T_j)=\lambda_i(T_j)$. Represent
  the ordering of
  the eigenvalues of $G$ by a bijection
  \begin{equation}
    \label{eq:sigma}
 \sigma\colon \{1,\dots,n\}\to H\text{ so that }d_{\sigma(p)}^*\ge
 d_{\sigma(q)}^*\text{ for
  }p\le q\text{ and so that }i\le q\text{ for }(j,i)=\sigma(q)
  \end{equation}
 (this is always
  possible, because the $d^*_i(T_j)$ are sorted
  nonincreasingly). Consider a diagram $D_G$ having $d_{\sigma(i)}^*$
  boxes in column $i$, then $\sum_{i=1}^k\lambda_i(G)$ counts the
  boxes in columns one to $k$. The Ferrers diagram $D_T$ of the
  intended threshold graph $T$ will be obtained
  from $D_G$ by only moving boxes to columns with smaller or equal
  index, then $\sum_{i=1}^k\lambda_i(T)=\sum_{i=1}^kd_i^*(T)\ge
  \sum_{i=1}^kd_{\sigma(i)}^*=\sum_{i=1}^k\lambda_i(G)$.

  For column $q=1,\dots,n$ in $D_G$ and $(j,i)=\sigma(q)$ place the $d_i^*(T_j)$
  boxes of this column by the following algorithm in the new diagram $D_T$. The box of row
  $r\in\{1,\dots,\min\{i,d_i^*(T_j)\}\}$ is
  placed on or above the diagonal of $D_T$, concretely in row $r$ in
  the next free column $c=q-|\{p\in\{1,\dots,q-1\}\colon (\bar\jmath,\bar\imath)=\sigma(p)
   \wedge \bar\imath<r\}|$, thus
  \eqref{eq:sigma} implies $r\le c\le q$; the box of row
  $r\in\{i+1,\dots,d_i^*(T_j)\}$ is placed in column $i$ (recall that
  $i\le q$ by \eqref{eq:sigma}) in the next free row $i+\sum_{(\hat
    \jmath,i)\in\{(\bar\jmath,\bar\imath)=\sigma(p)\colon
    \bar\imath=i \wedge p\in\{1,\dots,q-1\}\}}
  \max\{0,d_{i}^*(T_j)-i\}+(r-i)$, thus below the diagonal of
  $D_T$.

  We complete the proof by showing that $D_T$ is the Ferrers
  diagram of a threshold graph with trace $f=\max\{i\colon \exists
  (j,i)\in H \text{ with } d_i^*(T_j)>i\}$. Note that no boxes are
  placed on or above the diagonal of $D_T$ for rows $r>f$ (indeed,
$d_i^*(T_j)\neq i$ for all threshold graphs $T_j$ and $i$ due to the
transposed structure of their diagrams) and no boxes are placed
 below the diagonal for columns $c>f$.
Column $c=1,\dots,f$ contains
$\sum_{j=1}^h\max\{0,d_{c}^*(T_j)-c\}$ boxes below the
diagonal. The number of boxes in row $r=1,\dots,f$ of $D_T$ on or
above the diagonal computes to
$\sum_{j=1}^{h}|\{i\in\{r,\dots,n_j\}\colon d_i^*(T_j)\ge
r\}|=\sum_{j=1}^{h}\max\{0,d_r(T_j)-r+1\}=
\sum_{j=1}^{h}\max\{0,d_r^*(T_j)-r\}$, where the last equation
uses the defining property $d^*_i(T_j)=d_i(T_j)+1$ for
$i=1,\dots,f(T_j)$ (while $\max\{d^*_i(T_j), d_i(T_j)\}\le f(T_j)$ for
$i>f(T_j)$). Thus, for $c=r\in\{1,\dots,f\}$ the counts coincide
and $D_T$ is the Ferrers diagram of a threshold graph $T$ with $\sum_{i=1}^k\lambda_i(T)\ge\sum_{i=1}^k\lambda_i(G)$.
\end{Proof}
Forming the complement does not pose a problem as we show next.
\begin{Lemma}\label{L:complement}
Let $G$ be a graph on $n$ nodes with $m$ edges and suppose that for each $k\in\{1,\dots,n-1\}$ there is a threshold graph $T$ on $n$ nodes and $m$ edges
so that $\sum_{i=1}^k\lambda_i(T)\ge\sum_{i=1}^k\lambda_i(G).$
Then for the complement graphs $\bar G$ and $\bar T$ there holds
$\sum_{i=1}^{n-k-1}\lambda_i(\bar T)\ge\sum_{i=1}^{n-k-1}\lambda_i(\bar G)$. Because of
$\lambda_n=0$ any threshold graph $T$ having the same number of nodes
and edges satisfies $\sum_{i=1}^{n-1}\lambda_i(\bar T)=\sum_{i=1}^{n-1}\lambda_i(\bar G)$.
\end{Lemma}
\begin{Proof}
The equation $L(\bar G)=nI_n-\mathbf{1}\mathbf{1}^\top-L(G)$ gives rise to the well known relation $\lambda_i(\bar G)=n-\lambda_{n-i}(G)$ for $i=1,\dots,n-1$ (and $\lambda_n(\bar G)=0$ as usual). Because the sum of the eigenvalues of the Laplacian equals twice the number of edges of the graph, we get first
$\sum_{i=k+1}^{n-1}\lambda_{i}(T)= \sum_{i=1}^{n-k-1}\lambda_{n-i}(T)\le\sum_{i=1}^{n-k-1}\lambda_{n-i}(G)$ and then\\ $\displaystyle\sum_{i=1}^{n-k-1}\lambda_{i}(\bar T) = (n-k-1)n-\sum_{i=1}^{n-k-1} \lambda_{n-i}(T) \ge (n-k-1)n- \sum_{i=1}^{n-k-1} \lambda_{n-i}(G) = \sum_{i=1}^{n-k-1}\lambda_{i}(\bar G)$.
\end{Proof}
The two previous lemmas establish \cref{Th:preserve}.

Cographs are defined recursively as (i) $K_1$ is a cograph, (ii) the disjoint union of cographs is a cograph and (iii) the complement of a cograph is a cograph. Since $K_1$ is a threshold graph, in view of Lemmas \ref{L:union} and \ref{L:complement} we have proved \cref{Th:cographs},

\section{Equivalence with Brouwer's conjecture}\label{S:brouwer}

In this section we prove that, together with the Grone-Merris-Bai
theorem, Brouwer's conjecture is equivalent to conjecturing that every
graph is spectrally threshold dominated. Since threshold graphs
are known to satisfy Brouwer's conjecture, it is clear that any spectrally
threshold dominated graph $G$ satisfies Brouwer's
conjecture. The core of the proof is therefore to construct for
arbitrary $n$ and $m\le {n \choose 2}$ a threshold graph that attains
Brouwer's eigenvalue bound.

\begin{Proof}[of \cref{Th:brouwer}]
Note that by the Grone-Merris-Bai theorem Brouwer's conjecture is equivalent to $\sum_{i=1}^k\lambda_i(G)\le \min\{kn,m+k(k+1)/2,2m\}$ holding for $k\in\{1,\dots,n\}$, because no conjugate degree exceeds $n$ and the sum of all eigenvalues is $2m$.

Thus the equivalence is proven if for arbitrary $k\in\{1,\dots,n\}$ we show $\min\{kn,m+k(k+1)/2,2m\}=\max\{\sum_{i=1}^kd_i^*(T)\colon T\text{ threshold graph on $n$ nodes and $m$ edges}\}$. Depending on the relation between $k$, $n$ and $m$, we discern the following cases:

{\bf Case 1.} \underline{$\min\{kn,m+k(k+1)/2,2m\}=kn$:}
Consider the threshold graph $T$ constructed by filling up the Ferrers
diagram below the diagonal in columnwise order (on and above the
diagonal in corresponding rowwise order). The first $k$ columns below
the diagonal are fully filled because they require $kn-k(k+1)/2\le m$
boxes. Hence $T$  satisfies $d_i^*(T)=n$ for $i=\{1,\dots,k\}$ and
$\sum_{i=1}^k\lambda_i(T)=\sum_{i=1}^kd_i^*(T)=kn$. This is the
maximum attainable over all threshold graphs on $n$ nodes.

{\bf Case 2.} \underline{$\min\{kn,m+k(k+1)/2,2m\}=m+k(k+1)/2$:} In this case put $h:=\lfloor \tfrac{m}{k}+\tfrac{k+1}2\rfloor<n$ and $r:=m+k(k+1)/2-kh<k$. Note that this implies $h\ge k+1$. Define a threshold graph $T$ on $n$ nodes with $m$ edges of trace $k$ by the conjugate degrees
\begin{displaymath}
d_i^*(T)=\left\{\begin{array}{ll}
h+1 & i\le r,\\
h & r< i\le k,
\end{array}\right.
\end{displaymath}
then $\sum_{i=1}^k\lambda_i(T)=\sum_{i=1}^kd_i^*(T)=m+k(k+1)/2$. This value cannot be exceeded by any threshold graph on $n$ nodes with $m$ edges by the Grone-Merris-Bai Majorization theorem, because in the Ferrers diagram of the conjugate degrees up to column $k$ all boxes are used on and above the diagonal, while all possible $m$ boxes are included below the diagonal.

{\bf Case 3.} \underline{$\min\{kn,m+k(k+1)/2,2m\}=2m$:} Put $h:=\max\{h\in\{1,\dots,n\}:h(h+1)\le 2m\}<k$ and $r:=(2m-h(h+1))/2<h+1$, then the threshold graph $T$ of
trace $h$ with conjugate degrees
\begin{displaymath}
d_i^*(T)=\left\{\begin{array}{ll}
h+2 & i\le r,\\
h+1 & r< i\le h,\\
r& i=h+1,\\
0& h+1<i,
\end{array}\right.
\end{displaymath}
satisfies $\sum_{i=1}^k\lambda_i(T)=\sum_{i=1}^kd_i^*(T)=2m$ and this is the maximum attainable over all threshold graphs with $m$ edges.
\end{Proof}

\begin{figure}[h!]\begin{center}
\begin{tikzpicture}
[fill=black!150,scale=0.7,auto=left,every node/.style={circle,scale=0.7}]
  \foreach \i in {1,...,8}{
    \node[draw,circle] (\i) at ({360/8 * (\i - 1)+90}:2) {};} % :2 is the radius
\draw (1) -- (2) -- (3) -- (4) -- (7);
\draw (1) -- (3) -- (5) -- (7);
\draw (1) -- (4);
\draw (1) -- (5);
\draw (1) -- (6) -- (3);
\draw (2) -- (5) -- (7);
\draw (2) -- (4) -- (8);
\draw (2) -- (6);
  \foreach \i in {0,1/2,1,3/2,2}{
    \node[draw,rectangle] (\i) at (\i+4,2){};}
\foreach \i in {0,1/2,1,3/2,2}{
    \node[draw,rectangle] (\i) at (\i+4,1.5){};}
\foreach \i in {0,1/2,1,3/2,2}{
    \node[draw,rectangle] (\i) at (\i+4,1){};}
\foreach \i in {0,1/2,1,3/2,2}{
    \node[draw,rectangle] (\i) at (\i+4,.5){};}
\foreach \i in {0,1/2,1,3/2}{
    \node[draw,rectangle,fill] (\i) at (\i+4,2-\i) {};}
\foreach \i in {0,1/2,1,3/2}{
    \node[draw,rectangle] (\i) at (\i+4,0) {};}
\foreach \i in {0,1/2,1}{
    \node[draw,rectangle] (\i) at (\i+4,-.5) {};}
\foreach \i in {0,1/2}{
    \node[draw,rectangle] (\i) at (\i+4,-1){};}
\node[draw,rectangle] () at (4,-1.5){};
\end{tikzpicture}
\end{center}
\caption{Graph for illustration of the equivalence} \label{BroEx}
\end{figure}
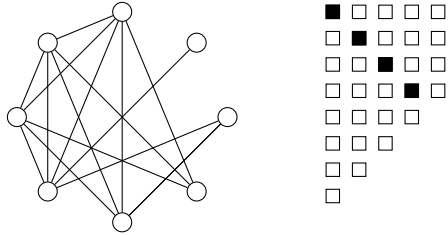
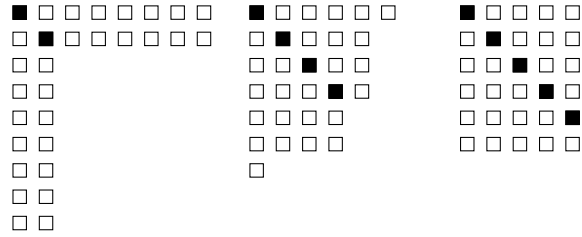
\begin{figure}[h!]\begin{center}
\begin{tikzpicture}
[fill=black!200,scale=0.7,auto=left,every node/.style={rectangle,scale=0.7}]
  \foreach \i in {0,1/2,1,3/2,2,5/2,3,7/2,4}{
    \node[draw,rectangle] (\i) at (0,-\i){};}
\foreach \i in {0,1/2,1,3/2,2,5/2,3,7/2,4}{
    \node[draw,rectangle] (\i) at (-.5,-\i){};}
\foreach \i in {0,1/2,1,3/2,2,5/2}{
    \node[draw] (\i) at (\i+.5,0){};}
\foreach \i in {1/2,1,3/2,2,5/2,3}{
    \node[draw] (\i) at (\i,-.5){};}
\node[draw,fill] () at (0,-.5){};
\node[draw,fill] () at (-.5,0){};

\foreach \i in {0,1/2,1,3/2,2,5/2,3}{
    \node[draw,rectangle] (\i) at (4,-\i){};}
\foreach \i in {0,1/2,1,3/2,2,5/2}{
    \node[draw,rectangle] (\i) at (4.5,-\i){};}
\foreach \i in {0,1/2,1,3/2,2,5/2}{
    \node[draw,rectangle] (\i) at (5,-\i){};}
\foreach \i in {0,1/2,1,3/2,2,5/2}{
    \node[draw,rectangle] (\i) at (5.5,-\i){};}
\foreach \i in {0,1/2,1,3/2}{
    \node[draw,rectangle] (\i) at (6,-\i){};}
\node[draw] () at (6.5,0) {};
\foreach \i in {0,1/2,1,3/2}{
    \node[draw,fill] (\i) at (4+\i,-\i){};}

\foreach \i in {0,1/2,1,3/2,2,5/2}{
    \node[draw,rectangle] (\i) at (8,-\i){};}
\foreach \i in {0,1/2,1,3/2,2,5/2}{
    \node[draw,rectangle] (\i) at (8.5,-\i){};}
\foreach \i in {0,1/2,1,3/2,2,5/2}{
    \node[draw,rectangle] (\i) at (9,-\i){};}
\foreach \i in {0,1/2,1,3/2,2,5/2}{
    \node[draw,rectangle] (\i) at (9.5,-\i){};}
\foreach \i in {0,1/2,1,3/2,2,5/2}{
    \node[draw] (\i) at (10,-\i){};}
\foreach \i in {0,1/2,1,3/2,2}{
    \node[draw,fill] (\i) at (8+\i,-\i){};}
\end{tikzpicture}

\end{center}
\caption{Transforms for min$\{kn,m+k(k+1)/2,2m\} = kn,m+k(k+1)/2,2m$}\label{transf}
\end{figure}

\begin{Example} Consider the graph and its Ferrers diagram of Figure \ref{BroEx}. There are $n=8$ vertices and $m=15$ edges. For $k=1,2$ we are in Case 1, for $k=3,4,5$ it is Case 2 and for $k=6,7,8$ we are in Case 3 of the theorem. We illustrate the construction of the threshold graphs $T$ for which $ \min \{kn,m+k(k+1)/2,2m\} = \max\{\sum_{i=1}^kd_i^*(T)\colon T\text{ threshold graph on $n$ nodes and $m$ edges}\}$ in Figure \ref{transf} for $k=2$, representing Case 1 (left), for $k=4$, representing Case 2 with $h=6,r=1$ (center) and for $k=7$, representing Case 3 with $h=5,r=0$ (right).
\end{Example}

We remark that in view of this result, whenever a class of graphs is shown to satisfy Brouwer's conjecture, the class also is spectrally threshold dominated and, in particular, the Laplacian energy is bounded by the  Laplacian energy of threshold graphs.

\begin{Corollary} Trees, unicyclic and bicyclic graphs are spectrally threshold dominated.
\end{Corollary}

This follows from the fact that in \cite{Haemers20102214}, it is proven that
trees satisfy Brouwer's conjecture. Likewise, in \cite{Du20123672}, it is proven that unicyclic and bicyclic graphs satisfy Brouwer's conjecture. Hence the Laplacian energy of these classes of graphs are also bounded by the Laplacian energy of threshold graphs. An explicit construction of a threshold graph on $n$ nodes and $m$ edges maximizing the Laplacian energy over all such threshold graphs is given in \cite{Helm15}. In the case of trees, it is known that the star on $n$ vertices (a threshold graph) has largest Laplacian energy among all trees with $n$ vertices \cite{Fri11}

In the same direction in \cite{Rocha201495} and in \cite{Wang201260}, it is proven that the conjecture of Brouwer holds for further classes of graphs.

\section*{Acknowledgments}

This work is partially supported by CAPES Grant PROBRAL 408/13 -
Brazil and DAAD PROBRAL Grant 56267227 - Germany. Trevisan also
acknowledges the support of CNPq - Grants 305583/2012-3 and
481551/2012-3.

%\bibliographystyle{abbrv}
%\bibliography{laplace}
%\end{document}

\end{document}